\documentclass[12pt]{article}

\newtheorem{theorem}{Theorem}[section]
\newtheorem{lemma}{Lemma}[section]
\newtheorem{definition}{Definition}[section]

\newtheorem{corollary}{Corollary}[section]
\newtheorem{proposition}{Proposition}[section]
\usepackage{setspace}
\usepackage{cite}
\usepackage{mathrsfs}
\usepackage{amsmath}
\usepackage{amssymb}
\usepackage{amsfonts}
\usepackage{amsmath,amsfonts,amssymb}
\usepackage{graphicx}
\usepackage{subfigure}
\usepackage{caption2}
\usepackage[rightcaption]{sidecap}
\usepackage{fancyhdr}
\newcommand{\R}{\mathbb R}

\newcommand{\B}{\mathbb B}
\newcommand{\N}{\mathbb N}

\newtheorem{remark}{Remark}[section]

\def \qed{\hfill{$\Box$}}

\begin{document}
\title{Some remarks on biharmonic elliptic problems with a singular nonlinearity
\thanks {Biharmonic operator,  singular nonlinearity, minimal
solutions, extremal solutions;}
\thanks { This work was supported by the the Natural Science Foundation of
China (Grant No: 10971061).}}
\author{Baishun Lai
}
\date{}
\maketitle
\begin{center}
\begin{minipage}{130mm}
{\small {\bf Abstract}\ \ \ We study the following semilinear
biharmonic equation
$$
\left\{
\begin{array}{lllllll}
\Delta^{2}u=\frac{\lambda}{1-u}, &\quad \mbox{in}\quad \B,\\
u=\frac{\partial u}{\partial n}=0, &\quad \mbox{on}\quad  \partial\B,\\
 \end{array}
\right.
$$
where $\B$ is the unit ball in $\R^{n}$ and $n$ is the exterior unit
normal vector. We prove the existence of $\lambda^{*}>0$ such that
for $\lambda\in (0,\lambda^{*})$ there exists a minimal (classical)
solution $\underline{u}_{\lambda}$, which satisfies
$0<\underline{u}_{\lambda}<1$. In the extremal case
$\lambda=\lambda^{*}$, we prove the existence of a weak solution
which is unique solution even in a very weak sense. Besides, several
new difficulties arise and many problems still remain to be solved.
we list those of particular interest in the final section. }
\end{minipage}
\end{center}
\vskip 0.2in \setcounter{equation}{0}
 \setcounter{section}{0}
\section{Introduction and results} 

\quad \quad In the last forty years a great deal has been written
about existence and multiplicity of solutions to nonlinear second
order elliptic problems in bounded and unbounded domains of $\R^{n}
(n\geq2)$. Important achievements on this topic have been made by
applying various combinations of analytical techniques, which
include the variational and topological methods. For the latter, the
fundamental tool which has been widely used is the maximum
principle. However, for higher order problems, a possible failure of
the maximum principle causes several technical difficulties, which
attracted the interest of many researchers. In particular, recently
fourth order equations with an singular non-linearity have been
studied extensively.  The motivation for considering these equations
stems from a model for the steady states of a simple micro
electromechanical system (MEMS) which has the general form (see for
example \cite{Pe})
$$
\left\{
\begin{array}{lllllll}
\alpha \Delta^{2}u=(\beta\int_{\Omega}|\nabla u|^{2} dx+\gamma)\Delta u+
\frac{\lambda f(x)}{(1-u)^{2}(1+\chi\int_{\Omega}\frac{dx}{(1-u)^{2}})}& \mbox{in}\ \ \Omega,\\
0<u<1&  \mbox{in}\ \ \Omega,\\
u=\alpha\frac{\partial u}{\partial n}=0 &  \mbox{on }\ \partial\Omega,
\end{array}
\right.
\eqno(M_{\lambda})
$$
where $\Delta^{2}(\cdot):=-\Delta(-\Delta)$ denotes the biharmonic
operator, $\Omega\subset\R^{n}$ is a smooth bounded domain, $n$
denotes the outward pointing unit normal to $\partial\Omega$ and
 $\alpha,\beta,\gamma,\chi\geq0,$ are physically relevant constants, $f\geq 0$ represents the
permittivity profile, $\lambda>0$ is a constant which is increasing
with respect to the applied voltage.\vskip 0.1in

Take $\alpha=\beta=\chi=0$ and $\gamma=1$, one obtain a simple
approximation of $(M_{\lambda})$
$$
\left\{
\begin{array}{lllllll}
-\Delta u=\lambda\frac{f(x)}{(1-u)^{2}} &\ \ \mbox{in}\ \Omega,\\
0<u<1  &\ \ \mbox{in}\ \Omega,\\
u=0   &\ \ \mbox{on}\ \partial\Omega.
\end{array}
\right.
\eqno(S_{\lambda})
$$
This simple model, which lends itself to the vast literature on
second order semilinear eigenvalue problems, is already a rich
source of interesting mathematical problems, see e.g.
\cite{Es,Es1,Gh} and the  references cited therein.

The case where $\gamma=\beta=\chi=0$ and $\alpha=1, f(x)\equiv 1$ in
the above model, that is when we replace $(1-u)^{-2}$ with
$(1-u)^{-p}$
$$
\left\{
\begin{array}{lllllll}
\Delta^{2} u=\frac{\lambda}{(1-u)^{p}} &\ \ \mbox{in}\ \Omega,\\
0<u<1  &\ \ \mbox{in}\ \Omega,\\
u=\frac{\partial u}{\partial n}=0   &\ \ \mbox{on}\ \partial\Omega.
\end{array}
\right.
\eqno(P_{\lambda})
$$
Because of the lack of a \textquotedblleft maximum principle ",
which play such a crucial role in developing the theory for the
Laplacian,  for $\Delta^{2}$ with Dirichlet boundary condition in
general domains (i.e., $\Omega\neq \B$),  very little is known about
$(P_{\lambda})$. As far as we are aware, only a paper \cite{Ghe}
study this problem for general domains. However, if $p>1$ and the
$\Omega$ is a ball, $(P_{\lambda})$ has recently been studied
extensively, see e.g. \cite{Ca,Co1,D,Gu4,Li,Mo} and its references.
One of the reasons to study $(P_{\lambda})$ in a ball is that a
maximum principle holds in this situation, see \cite{Bo}, and so
some tools that are well suited for $(S_{\lambda})$ can work for
$(P_{\lambda})$. The second reason is that one can easily find a
explicit singular radial solution,  denoted by $1-
|x|^{\frac{4}{p+1}} (p>1)$, of $(P_{\lambda})$ for $\Omega=\B$ and a
suitable parameter $\lambda$ which satisfy the first boundary
condition but not the second. The singular radial solution, called
\textquotedblleft ghost" singular solution, play a fundamental role
to characterize the \textquotedblleft true" singular solution, see
in particular \cite{Co1}.

In this paper, we will focus essentially our attention on the case
where $p=1$ and $\Omega$ is a ball, namely
$$
\left\{
\begin{array}{lllllll}
\Delta^{2}u=\frac{\lambda}{1-u} & \mbox{in}\  \B,\\
0<u\leq 1 &  \mbox{in}\ \ \B, \\
u=\frac{\partial u}{\partial n} =0 &  \mbox{on}\  \  \partial \B.\\
\end{array}
\right.
\eqno(1.1)_{\lambda}
$$
For the corresponding second order problem, which is related to the
general study of singularities of minimal hypersurfaces of Euclidean
space, has been studied by Meadows, see \cite{Mea}.  In that case,
however, the start point was an explicit singular solution (i.e.,
$u_{s}(r)=r^{2}$) with parameter $\lambda=n-1$. When turning to the
biharmonic problem $(1.1)_{\lambda}$, one can not find any explicit
singular solution even \textquotedblleft ghost" singular solution
which causes several technical difficulties.  The first purpose of
the present paper is to extend $(1.1)_{\lambda}$ some well-known
results relative to $(P_{\lambda}) $. The second (and perhaps most
important) purpose of the present paper is to emphasize some
striking differences between $(1.1)_{\lambda}$ and $(P_{\lambda}) $.
\vskip0.2in

\subsection{ Preliminaries}

\quad \quad Besides classical solution i.e. $u\in C^{4}(\bar{\B})$
which satisfy $(1.1)_{\lambda}$, let us introduce the class of weak
solutions we will be dealing with. We denote by $H_{0}^{2}(\B)$ the
usual Sobolev space which can be defined by completion as follows:
$$
H_{0}^{2}(\B):=cl\{u\in C_{0}^{\infty}(\B): \|\Delta u\|_{2}<\infty\}
$$
and which is an Hilbert space endowed with the scalar product
$$
(u,v)_{H_{0}^{2}(\B)}:=\int_{\B}\Delta u\Delta v dx
$$

\begin{definition}\label{D1.1}
 We say that $u\in L^{1}(\B)$ is
a weak solution of $(1.1)_{\lambda}$ provided $0\leq u\leq 1$ almost
everywhere, $\frac{1}{1-u}\in L^{1}(\B)$ and
\begin{equation}\label{Equ1}
\int_{\B} u\Delta^{2}\varphi dx=\lambda \int_{\B}\frac{\varphi}{(1-u)}dx,\ \ \forall \varphi\in C^{4}(\bar{\B})\cap H_{0}^{2}(\B) 
\end{equation}
\end{definition}

When in (\ref{Equ1}) the equality is replaced by the inequality
$\geq$ (resp.$\leq$) and $\varphi\geq0$, we say that $u$ is a weak
super-solution (resp. weak sub-solution) of (\ref{Equ1}) provided
the following boundary conditions are satisfied: $u=0 $ (resp.$=$)
and $\frac{\partial u}{\partial n}\leq 0$ (resp.$\geq$) on $\partial
\B$. \vskip 0.1in

\begin{definition}\label{D1.1}
We call a solution $u$ of $(1.1)_{\lambda}$ minimal if $u\leq v$
a.e. in $\B$ for any further solution $v$ of $(1.1)_{\lambda}$
\end{definition}

If $u$ is a classical solution of $(1.1)_{\lambda}$, then it turns
out to be well defined the linearized operator at $u$
$$
L_{u}:=\Delta^{2}-\frac{\lambda}{(1-u)^{2}}
$$
which yields the following notion of stability \vskip 0.1in

\begin{definition}\label{D1.3}
 A classical solution $u$ of
$(1.1)_{\lambda}$ is semi-stable provided
$$
\mu_{1}(u)=\inf\left\{\int_{\B}(\Delta \varphi)^{2}-\frac{\lambda\varphi^{2}}{(1-u)^{2}}:
\phi\in H_{0}^{2}(\B), \|\phi\|_{L^{2}}=1\right\}\geq0
$$
If $\mu_{1}(u)>0$ we say that $u$ is stable.
\end{definition}


As far as we are concerned with weak solutions, the linearized
operator is no longer well defined,
 however we introduce the following weaker notion of stability.\vskip 0.1in

\begin{definition}\label{D1.4}
 A weak solution $u$ to $(1.1)_{\lambda}$ is said to be weakly
stable if $\frac{1}{(1-u)^{2}}\in L^{1}(\B)$ and the following
holds:
$$
\int_{\B}|\Delta\varphi|^{2}dx\geq \int_{\B}\frac{\lambda\varphi^{2}}{(1-u)^{2}}, \varphi\in H_{0}^{2}(\B), \varphi\geq0.
$$
\end{definition}

According to the class of solutions which we consider, let us
introduce the following values:
\begin{equation}\label{Equ2}
\begin{array}{lllllll}
\lambda^{*}:=\sup\{\lambda\geq0: (1.1)_{\lambda}\ \ \mbox{ posses a weak solution}\};\\
\lambda_{*}:=\sup\{\lambda\geq0: (1.1)_{\lambda}\ \ \mbox{ posses a classical solution}\}.\\
\end{array}
\end{equation}

\begin{remark}\label{R1.1}
 Clearly, a classical solution is also
a weak solution, so that one has $\lambda_{*}\leq \lambda^{*}$.
Moreover, by standard elliptic regularity theory for the biharmonic
operator \cite{Ag}, any weak solution of $(1.1)_{\lambda}$ which
satisfies $\|u_{\lambda}\|< 1$ turns out to be smooth.
\end{remark}


Besides, we give a notion of $H_{0}^{2}$($\B$)- weak solutions,
which is an intermediate class between classical and weak solutions.

\begin{definition}\label{D1.5}
 We say that $u$ is a
$H_{0}^{2}$($\B$)- weak solution of (1.1) if  $(1-u)^{-1}\in
L^{1}(\B)$ and if
$$
\int_{\B}\Delta u\Delta\phi=\lambda\int_{\B}\phi(1-u)^{-1},\ \ \ \forall \phi\in C^{4}(\bar{\B})\cap H_{0}^{2}(\B).
$$
 We say that $u$ is a $H_{0}^{2}$($\B$)- weak super-solution (resp. $H_{0}^{2}$($\B$)- weak sub-solution) of $(1.1)_{\lambda}$
  if for $\phi\geq0$ the equality is replaced with $\geq$ (resp.$\leq$) and $u\geq 0$ (resp. $\leq$), $\frac{\partial u}{\partial n}\leq 0$
  (resp. $\geq$) on $\partial \B$.
\end{definition}

\subsection{Main results}

\quad \quad In order to state our results, we denote by $\nu_{1}$
the first eigenvalue of the biharmonic operator on $\B$ with
Dirichlet boundary conditions, which is characterized variationally
as follows:
$$
\nu_{1}:=\inf\left\{\int_{\B}|\Delta u|^{2}dx:\ \ u\in H_{0}^{2}(\B), \|u\|=1\right\}.
$$
It is well known that $\nu_{1}>0$, that it is simple, isolated and
that the corresponding eigenfunctions $\psi>0$, spherically
symmetric , radially decreasing and do not change sign.

We may now state the following theorem. \vskip 0.1in

\begin{theorem}\label{result1}
  There exists $\lambda_{*}>0$ such that
for $0<\lambda<\lambda_{*}$, $(1.1)_{\lambda}$ poses a minimal
classical solution, denoted by $\underline{u}_{\lambda}$, which is
positive and stable. Moreover, $\lambda_{*}$ satisfies the following
bounds:
$$
\max\{4n(n-2), 2n(n+2)\}\leq \lambda_{*}\leq \frac{\nu_{1}}{4}.
$$
\end{theorem}

It is remarkable that at $\lambda_{*}$ there is an immediate switch
from existence of regular minimal solutions to nonexistence of any
(even singular) solution. The only possibly singular minimal
solution corresponds to $\lambda=\lambda_{*}$. This result is known
from \cite{Br} for the second order problem($S_{\lambda}$), but the
method used there may not be carried over to fourth order problems.
Nevertheless, the result extends to biharmonic case in the following
theorem.
\begin{theorem}\label{result2}
\  \  The following holds:
$$
\lambda_{*}=\lambda^{*}.
$$
In particular, for $\lambda>\lambda^{*}$ there are no solutions,
even in the weak sense. Furthermore, for almost every $x\in \B$,
there exists
$$
u^{*}(x):=\lim_{\lambda\to\lambda^{*}}u_{\lambda}(x)
$$
and $u^{*}(x)$ is a weakly stable $H_{0}^{2}(\B)-$ weak solution of
$(1.1)_{\lambda}$, which is called the extremal solution.

If $n\leq 4$ then the extremal solution $u^{*}$ of $(1.1)_{\lambda}$
is smooth, i.e.,
$u^{*}=\lim_{\lambda\to\lambda^{*}}\underline{u}_{\lambda}(x)$
exists in the topology of $C^{4}(\B)$. It is the unique regular
solutions to $(1.1)_{\lambda^{*}}$. \vskip 0.1in
\end{theorem}

 From the above theorem, we note that the function $u^{*}$ exists in any dimension, dose
solve $(1.1)_{\lambda^{*}}$ in the $H_{0}^{2}(\B)$ weak sense and it
is a classical solution in dimensions $1\leq n\leq 4$. This will
allow us to start another branch of nonminimal (unstable )
solutions. Besides, inspired by \cite{Ca,Ma,Da} we get the following
uniqueness of the extremal solution of $(1.1)_{\lambda}$, which
gives Theorem1.3.

\vskip 0.1in

\begin{theorem} \label{result3}
 \  \ Let $v$ be a weak super-solution of $(1.1)_{\lambda}$ with parameter
$\lambda^{*}$. Then $v=u^{*}$; in particular $(1.1)_{\lambda}$ has a
unique weak solution. \vskip 0.1in
\end{theorem}

From this theorem, we know that there are no strict super-solutions
to equation $(1.1)_{\lambda^{*}}$. \vskip 0.1in

\begin{corollary} \label{C1.1}
\ \ Let $u_{\lambda}\in H_{0}^{2}(\B)$ be a weak solution of
$(1.1)_{\lambda}$ such that $\|u_{\lambda}\|=1$. Then $u_{\lambda}$
is weakly stable if and only if $\lambda=\lambda^{*}$ and
$u_{\lambda}=u^{*}$
\end{corollary}

We may also characterize the uniform convergence to 0 of
$\underline{u}_{\lambda}$ as $\lambda\to 0$ by giving the precise
rate of its extinction.

\begin{theorem} \label{result4}
 For all $\lambda\in
(0,\lambda^{*})$ let $\underline{u}_{\lambda}$ be the minimal
solution of $(1.1)_{\lambda}$ and let
$$
V_{\lambda}(x)=\frac{\lambda}{8n(n+2)}\left[1-|x|^{2}\right]^{2}.
$$
Then $\underline{u}_{\lambda}>V_{\lambda}(x)$ for all
$\lambda<\lambda^{*}$ and all $|x|<1$, and
$$
\lim_{\lambda\to0}\frac{\underline{u}_{\lambda}}{V_{\lambda}(x)}=1\  \mbox{uniformly with respect to}\ \ x\in \B.
$$
\end{theorem}

\subsection{ Key-ingredients}

\quad \quad Now we give some comparison principles which will be
used throughout the paper

\begin{lemma}\label{L1.1}
\ (Boggio's principle, \cite{Bo}) If $u\in C^{4}(\bar{\B}_{R})$
satisfies
$$
\left\{
\begin{array}{lllllll}
\Delta^{2}u\geq 0 & \mbox{in}\ \  \B_{R},\\
u=\frac{\partial u}{\partial n}=0 &  \mbox{on}\ \  \partial \B_{R},
\end{array}
\right.
$$
then $u\geq 0$ in $\B_{R}$.
\end{lemma}

\begin{lemma}\label{L1.2}
Let $u\in L^{1}(\B_{R})$ and suppose that
$$
\int_{\B_{R}}u\Delta^{2}\varphi\geq0
$$
for all $\varphi\in C^{4}(\bar{\B}_{R})$ such that $\varphi\geq0$ in
$\B_{R}$, $\varphi|_{\partial \B_{R}}=\frac{\partial
\varphi}{\partial n}|_{\partial \B_{R}}=0$. Then $u\geq 0$ in
$\B_{R}$. Moreover $u\equiv 0$ or $u>0$ a.e., in $\B_{R}$.
\end{lemma}

\vskip0.1in

For a proof see Lemma 17 in \cite{Ar1}. From this lemma, we know
that any solution of $(1.1)_{\lambda}$ is necessarily positive a.e.
inside the ball.

\vskip0.1in \noindent
\begin{lemma}\label{L1.3}
\ If $u\in H^{2}(\B_{R})$ is
radial, $\Delta^{2}u\geq 0$ in $\B_{R}$ in the weak sense, that is
$$
\int_{\B_{R}}\Delta u\Delta\varphi\geq 0 \ \ \forall \varphi \in C_{0}^{\infty}(\B_{R}), \ \varphi\geq0
$$
and $u|_{\partial \B_{R}}\geq0, \frac{\partial u}{\partial
n}|_{\partial \B_{R}}\leq 0$ then $u\geq 0$ in $\B_{R}$.\vskip0.1in
\end{lemma}
{\bf Proof.} For the sake of completeness, we include a brief proof
here.  We only deal with the case $R=1$ for simplicity. Solve
$$
\left\{
\begin{array}{lllllll}
\Delta^{2}u_{1}=\Delta^{2}u & \mbox{in} \ \  \B\\
u_{1}=\frac{\partial u_{1}}{\partial n}=0 &  \mbox{on}\ \ \partial \B
\end{array}
\right.
$$
in the sense $u_{1}\in H_{0}^{2}(\B)$ and $\int_{\B}\Delta
u_{1}\Delta\varphi=\int_{\B}\Delta u\Delta\varphi$ for all $\varphi
\in C_{0}^{\infty}(\B)$. Then $u_{1}\geq 0$ in $\B$ by Lemma 2.2.

Let $u_{2}=u-u_{1}$ so that $\Delta^{2}u_{2}=0$ in $\B$. Define
$f=\Delta u_{2}$. Then $\Delta f=0$ in $\B$ and since $f$ is radial
we find that $f$ is a constant. It follows that $u_{2}=ar^{2}+b$.
Using the boundary conditions we deduce $a+b\geq 0$ and $a\leq0$,
which imply $u_{2}\geq0$.\vskip 0.1in

\begin{lemma} \label{L1.4}
Let $f\in L^{1}(\B_{R}), f\geq0$ almost everywhere. Then there
exists a unique $u\in L^{1}(\B_{R})$ such that $u\geq0$ and
\begin{equation}\label{Equ3}
\int_{\B_{R}}u \Delta^{2}\varphi=\int_{\B_{R}} f \varphi,\ \  \varphi\in C^{4}(\bar{\B}_{R})\cap H_{0}^{2}(\B_{R}). 
\end{equation}
Moreover, there exists $C>0$ which does not depend on $f$ such that
$\|u\|_{1}\leq C\|f\|_{1}$.  \vskip0.1in
\end{lemma}

{\bf Proof.}  The proof is standard, see \cite {Ar1}, we give a
proof here for the sake of completeness. The uniqueness is clear.
Indeed, let $v_{1}$ and $v_{2}$ be two solutions of (\ref{Equ3}).
Then $v=v_{1}-v_{2}$ satisfies
$$
\int_{\B}v\Delta\varphi=0  \ \  \varphi\in C^{4}(\bar{\B}_{R})\cap H_{0}^{2}(\B_{R}).
$$
Given any $\zeta \in C_{0}^{\infty}(\B)$ let $\varphi$ be the
solution of
$$
\left\{
\begin{array}{lllllll}
\Delta^{2}\varphi= \zeta & \mbox{in}\ \  \B,\\
\varphi=\frac{\partial \varphi}{\partial n}=0 &  \mbox{on}\ \  \partial \B.
\end{array}
\right.
$$
It follows that
$$
\int_{\B}v\zeta=0.
$$
Since $\zeta$ is arbitrary, we deduce that $v=0$.

For the existence, Given an integer $k\geq0$ we set
$f_{k}=\min\{f(x),k\}$, so that $f_{k}\rightarrow f $ as
$k\to\infty$ in $L^{1}(\B)$. Let $v_{k}$ be the solution of
\begin{equation}\label{Equ4}
\left\{
\begin{array}{lllllll}
\Delta^{2}v_{k}= f_{k} & \mbox{in}\ \  \B,\\
v_{k}=\frac{\partial v_{k}}{\partial n}=0 &  \mbox{on}\ \  \partial \B.
\end{array}
\right.
\end{equation}
The sequence $(v_{k})_{k\geq0}$ is clearly monotone nondecreasing.
It is also a cauchy sequence in $L^{1}(\B)$ since
$$
\int_{\B}(v_{k}-v_{l})=\int_{\B}(f_{k}-f_{l})\zeta_{0},
$$
where $\zeta_{0}$ is defined by
$$
\left\{
\begin{array}{lllllll}
\Delta^{2}\zeta_{0}= 1 & \mbox{in}\ \  \B,\\
\zeta_{0}=\frac{\partial \zeta_{0}}{\partial n}=0 &  \mbox{on}\ \  \partial \B.
\end{array}
\right.
$$
Hence
$$
\int_{\B}|v_{k}-v_{l}|\leq C\int_{\B}|f_{k}-f_{l}|dx.
$$
Passing to the limit in (\ref{Equ4}) (after multiplication by
$\varphi$) we obtain (\ref{Equ3}) and $u\geq0$ according to the
Lemma 1.2. Finally, taking $\varphi=\zeta_{0}$ in (\ref{Equ3}), we
obtain
$$
\|v\|_{L^{1}}=\int_{\B}v=\int_{\B}f\zeta_{0}\leq C\|f\|_{L^{1}}
$$
and the proof is completed.

\begin{proposition}\label{P1.1}
 Assume the existence of
a weak super-solution $U$ of $(1.1)_{\lambda}$. Then there exists a
weak solution $u$ of $(1.1)_{\lambda}$ so that $0\leq u\leq U$ a.e
in $\B$.
\end{proposition}

{\bf Proof.}\ \  By means of a standard monotone iteration argument,
set $u_{0}:= U$ and define recursively $u_{n+1}\in L^{1}(\B)$ as the
unique solution of
$$
\int_{\B}u_{n+1}\Delta^{2}\varphi dx=\lambda\int_{\B}\frac{\varphi}{(1-u_{n})^{2}}dx, \varphi\in C^{4}(\bar{\B})\cap H_{0}^{2}(\B),
$$
then we have
$$
\int_{\B}(u_{n}-u_{n+1})\Delta^{2}\varphi dx\geq 0, \varphi\in C^{4}(\bar{\B})\cap H_{0}^{2}(\B)
$$
and Lemma \ref{L2.2} yields $0\leq u_{n+1}\leq u_{n}<U(x)$ a.e. for
all $n\in \N$. Since
$$
(1-u_{n})^{-1}\leq (1-U)^{-1}\ \ \in L^{1}(\B),
$$
and the claim follows from the Lebesgue convergence Theorem.
\vskip0.1in

We complete these preliminary results by proving a key lemma which
provides a comparison principle.\vskip0.1in
\begin{lemma}\label{L1.5}
 \ Assume $u_{1}$ is a weakly stable
$H_{0}^{2}$($\B$)- weak
 sub-solution of $(1.1)_{\lambda}$ and $u_{2}$ is
$H_{0}^{2}$($\B$)- weak super-solution of $(1.1)_{\lambda}$. Then,

 (1) $u_{1}\leq u_{2}$ almost
everywhere in $\B$.

(2)  if $u$ is a classical solution such that $\mu_{1}(u)=0$ and $U$
is any classical super-solution of  $(1.1)_{\lambda}$, then $u\equiv
U$.
\end{lemma}

{\bf Proof.}\ \ (1)  Define $\omega:=u_{1}-u_{2}$. Then by the
Moreau decomposition \cite{Mo} for the biharmonic operator, there
exists $\omega_{1}, \omega_{2}\in H_{0}^{2}(\B)$, with
$\omega=\omega_{1}+\omega_{2}, \omega_{1}\geq0$ a.e.,
$\Delta^{2}\omega_{2}\leq0$ in the $H_{0}^{2}(\B)-$ weak sense and
$$
\int_{\B}\Delta \omega_{1}\Delta \omega_{2}=0
$$
By Lemma \ref{L1.1}, we have that $\omega_{2}\leq 0$ a.e. in $\B$.

Given now $0\leq\varphi\in C_{0}^{\infty}(\B)$, we have that
$$
\int_{\B}\Delta\omega\Delta\varphi\leq\lambda\int_{\B}(f(u_{1})-f(u_{2}))\varphi,
$$
where $f(u)=(1-u)^{-1}$. Since $u$ is stable, one has
$$
\lambda\int_{\B}f'(u)\omega_{1}^{2}\leq \lambda\int_{\B}(\Delta\omega_{1})^{2}=\lambda\int_{\B}\Delta\omega\Delta\omega_{1}
\leq\lambda\int_{\B}(f(u_{1})-f(u_{2}))\omega_{1}
$$
Since $\omega_{1}\geq \omega$ one also has
$$
\int_{\B}f'(u)\omega \omega_{1}\leq \int_{\B}(f(u_{1})-f(u_{2}))\omega_{1}
$$
which once re-arrange gives
$$
\int_{\B}\tilde{f}\omega_{1}\geq 0,
$$
where $\tilde{f}(u_{1})=f(u_{1})-f(u_{2})-f'(u_{1})(u_{1}-u_{2})$.
The strict convexity of $f$ gives $\tilde{f}\leq0$ and $\tilde{f}<0$
whenever $u\neq U$. Since $\omega_{1}\geq0$ a.e. in $\B$ one sees
that $\omega\leq0$ a.e. in $\B$. The inequality $u_{1}\leq u_{2}$
a.e. in $\B$ is then established.\vskip0.1in

(2)\  Let $\varphi>0$ be the first eigenfunction of
$\Delta^{2}-\lambda f'(u)$ in $H_{0}^{2}(\B)$, we now, for $0\leq
t\leq 1$, define
$$
g(t)=\int_{\B}\Delta(tU+(1-t)u)\Delta\phi-\lambda\int_{\B}f(tU+(1-t)u)\phi,
$$
where $\phi$ is the above first eigenfunction. Since $f$ is convex
one sees that
$$
g(t)\geq\lambda \int_{\B}[tf(U)+(1-t)f(u)-f(tU+(1-t)u)]\phi\geq0
$$
for every $t\geq0$. Since $g(0)=0$ and
$$
g'(0)=\int_{\B}\Delta(U-u)\Delta\phi-\lambda f'(u)(U-u)\phi=0,
$$
we get that
$$
g''(0)=-\lambda\int_{\B}f''(u)(U-u)^{2}\phi\geq0.
$$
Since $f''(u)\phi> 0$ in $\B$, we finally get that $U=u$ a.e. in
$\B$.\qed \vskip0.1in


\setcounter{equation}{0}
 \section{Existence results: proofs of Theorem \ref{result1} and
\ref{result2}}

\vskip0.1in

 \subsection{ The branch of minimal solutions}

 \quad \quad Let us define
 $$
 \Lambda:=\{\lambda\geq0: (1.1)_{\lambda}\ \mbox{has a classical solution with parameter}\ \lambda \}.
 $$

\begin{proposition}\label{P2.1}
 For all
$0\leq\lambda<\lambda_{*}$, there exists a minimal classical
solution $\underline{u}_{{\lambda}}$ of $(1.1)_{\lambda}$ which is
smooth and stable. Moreover, \vskip0.1in

(i) The map $\lambda\rightarrow \underline{u}_{{\lambda}},$ for
$\lambda\in (0,\lambda_{*})$ is differentiable and strictly
increasing;

(ii) The map $\lambda\rightarrow \mu_{1}(\underline{u}_{{\lambda}})$
is decreasing on $(0,\lambda_{*})$;

(iii) Let $\tilde{u}_{\lambda}$ be a regular solution of
$(1.1)_{\lambda}$ for $\lambda\in (0,\lambda_{*})$, if
$\tilde{u}_{\lambda}$ is not the minimal solution, then
$\mu_{1}(\tilde{u}_{\lambda})<0$.
\end{proposition}

 {\bf Proof.}\ \  First we show that $\Lambda$ dose not
consist of just $\lambda=0$. To this end, let $\psi_{R}$ be the
first eigenfunction of the biharmonic operator subject to Dirichlet
boundary conditions on $\B_{R}\supset \B$ which we normalize by
$\sup_{\B_{R}}\Psi_{R}=1$ and let $\nu_{R}>0$ be the corresponding
eigenvalue. Next, we are going to prove that for $\theta\in (0,1)$
the function $\psi=\theta\psi_{R}$ is a super-solution of
$(1.1)_{\lambda}$ as long as $\lambda$ is sufficiently small. We
have
$$
0<1-\theta\psi_{R}<1, \ \ \mbox{in}\ \ \B
$$
and moreover
$$
\Delta^{2}\psi=\nu_{R}\theta\psi_{R}\geq \frac{\lambda}{1-\theta\psi_{R}}=\frac{\lambda}{1-\psi}
$$
provide that
$$
\nu_{R}\theta\psi_{R}(1-\theta\psi_{R})\geq \lambda.
$$

Notice that
$$
0<s_{1}:=\inf_{x\in\B}\psi<s_{2}:=\sup_{x\in\B}\psi<1
$$
and that $\frac{\partial \psi}{\partial n}<0$ on $\partial\B$. Thus,
looking at the function $g(s)=s(1-s)$, for $s\in [s_{1}, s_{2}]$, it
is easily seen that we can choose $\lambda>0$ sufficiently small
such that
$$
\nu_{R}\inf_{x\in\B}g(\theta\psi(x))>\lambda.
$$
Since $\underline{u}\equiv0$ is a sub-solution of $(1.1)_{\lambda}$,
the classical sub-super solution Theorem provides a classical
solution $u_{\lambda}$ to $(1.1)_{\lambda}$. With such function
$u_{\lambda}$, we can use the Boggio principle to show
straightforwardly that the iterative scheme
\begin{equation}\label{Equ5}
\left\{
\begin{array}{lllllll}
\Delta^{2} u_{n,\lambda} =\frac{\lambda}{(1-u_{n,\lambda})}& \ \ \mbox{in}\ \B,   \\\\
 u_{n,\lambda}=\frac{\partial u_{n,\lambda}}{\partial n}  =0 & \ \ \mbox{in}\ \partial\B,  \\\\
 u_{0,\lambda}=0     & \ \ \mbox{in}\ \B,
 \end{array}
\right.
\end{equation}
gives rise to a monotone sequence $\{u_{n,\lambda}\}$ satisfying
$$
0=u_{0,\lambda}\leq u_{1,\lambda}  ... u_{n-1,\lambda}\leq. . .\leq u_{\lambda}<1
$$
for all $n\in \N$. Therefore the minimal solution
$\underline{u}_{\lambda}$ is obtained as the increasing limit
$$
\underline{u}_{\lambda}(x):=\lim_{n\rightarrow\infty}u_{n,\lambda}
$$
Again from the Boggio positivity preserving property (Lemma 1.1) we
obtain $0<\underline{u}_{\lambda}<1$; in particular, from standard
elliptic regularity theory for the biharmonic operator follows that
$\underline{u}_{\lambda}(x)$ is smooth. In order to prove stability,
let us argue as follows: set
$$
\lambda_{**}:=\sup\{\lambda\in (0,\lambda_{*}): \mu_{1}(\underline{u}_{\lambda})>0\}
$$
clearly $\lambda_{**}\leq\lambda_{*}$. Now suppose by contradiction
that $\lambda_{**}<\lambda_{*}$ and let $\varepsilon>0$ sufficiently
small such that $\lambda_{**}+\varepsilon<\lambda_{*}$ and
$v_{\lambda_{**}+\varepsilon}$ be the corresponding minimal
solution. By the definition and left continuity of the map
$\lambda\to \mu_{1}(\underline{u}_{\lambda})$ we have necessarily
$\mu_{1}(\underline{u}_{\lambda_{**}})=0$. Since
$v_{\lambda_{**}+\varepsilon}$ is a super-solution of
$(1.1)_{\lambda_{**}}$, by Lemma \ref{L1.5} we get
$v_{\lambda_{**}+\varepsilon}=\underline{u}_{\lambda_{**}}$ and thus
$\varepsilon=0$, a contradiction.

Since each $\underline{u}_{\lambda}$ is stable, then by setting
$F(\underline{u}_{\lambda},\lambda):=-\Delta^{2}-\frac{\lambda}{1-\underline{u}_{\lambda}}$,
we get that
$F_{\underline{u}_{\lambda}}(\underline{u}_{\lambda},\lambda)$ is
invertible for $0<\lambda<\lambda_{*}$. It then follows from
Implicit Function Theorem that $\underline{u}_{\lambda}(x)$ is
differentiable with respect to $\lambda$.

Now we prove the map $\lambda\to \underline{u}_{\lambda}$ is
strictly increasing on $(0,\lambda_{*})$. Consider
$\lambda_{1}<\lambda_{2}<\lambda_{*}$, their corresponding minimal
positive solutions $\underline{u}_{\lambda_{1}}$ and
$\underline{u}_{\lambda_{2}}$, and let $u^{*}$ be a solution for
$(1.1)_{\lambda_{2}}$. The same as the above iterative scheme, we
have
$$
\underline{u}_{\lambda_{1}}=\lim_{n\to\infty}u_{n}(\lambda_{1};x)\leq u^{*}\ \ \mbox{in}\ \ \B,
$$
and in particular
$\underline{u}_{\lambda_{1}}\leq\underline{u}_{\lambda_{2}}$ in
$\B$. Therefore, $\frac{d\underline{u}_{\lambda}}{d\lambda}\geq0$
for all $x\in \B$.

Finally, by differentiating $(1.1)_{\lambda}$  with respect to
$\lambda$, and since $\lambda\to \underline{u}_{\lambda}$ is
nondecreasing, we get
$$
-\Delta^{2}\frac{d\underline{u}_{\lambda}}{d\lambda}-\frac{\lambda}{(1-\underline{u}_{\lambda})^{2}}\frac{d\underline{u}_{\lambda}}{d\lambda}
=\frac{\lambda}{1-\underline{u}_{\lambda}}\geq 0, \ x\in \B; \ \ \frac{d\underline{u}_{\lambda}}{d\lambda}=0, \ \ x\in \partial\B.
$$
Applying the strong maximum principle, we conclude that
$\frac{d\underline{u}_{\lambda}}{d\lambda}>0$ on $\B$ for all
$0<\lambda<\lambda_{*}$

That $\lambda\to \mu_{1,\lambda}$ is decreasing follow easily from
the variational characterization of $\mu_{1,\lambda}$, the
monotonicity of $\lambda\to \underline{u}_{\lambda}$, as well as the
monotonicity of $(1-\underline{u}_{\lambda})^{-2}$ with respect to
$\underline{u}_{\lambda}$ and the proof of the (ii) is completed.

Now we give the proof of (iii). Let $\underline{u}_{\lambda}$ be the
minimal solution for $(1.1)_{\lambda}$ so that
$\tilde{u}_{\lambda}\geq \underline{u}_{\lambda}$. If the
linearization around $\tilde{u}_{\lambda}$ had nonnegative first
eigenvalue, then Lemma \ref{L1.5} would also yield
$\tilde{u}_{\lambda}\leq \underline{u}_{\lambda}$ so that
$\tilde{u}_{\lambda}$ and $\underline{u}_{\lambda}$ necessarily
coincide, a contradiction.

\begin{remark}
Dose (iii) of Proposition (\ref{P2.1}) extend to weak solutions $u$
as formulated in [21, Theorem 3.1]?
\end{remark}

\vskip0.1in

\noindent
  \subsection{ Weak solutions versus classical solutions}

\begin{lemma}\label{L2.1}
  Let $u_{\mu}$ be a weak
solution of $(1.1)_{\mu}$ with $\mu<\lambda^{*}$. Then, for
$\varepsilon>0$ sufficiently small, the problem
$(1.1)_{(1-\varepsilon)\mu}$ posses a classical solution.
\end{lemma}
{\bf Proof.} Let $\tilde{u}\in L^{1}(\B)$ be the unique solution of
$$
\int_{\B}\tilde{u}\Delta^{2}\varphi=\mu\int_{\B}\frac{(1-\varepsilon)}{1-u_{\mu}}\varphi dx,\ \
\varphi\in C^{4}(\bar{\B})\cap H_{0}^{2}(\B)
$$
provided by Lemma \ref{L1.4}. By hypothesis we have
$$
\int_{\B}u_{\mu}\Delta^{2}\varphi dx=\mu\int_{\B}\frac{\varepsilon}{1-u_{\mu}}\varphi dx, \ \
\varphi\in C^{4}(\bar{\B})\cap H_{0}^{2}(\B).
$$
By uniqueness we get
$$
(1-\varepsilon)u_{\mu}=\tilde{u}
$$
whereas Lemma \ref{L1.2} yields $\tilde{u}>0$ a.e. in $\B$ and hence
we may assume
$$
u_{\mu}>\tilde{u}, \ \ x\in\B\setminus\{x\in\B: \tilde{u}=0\}
$$
Therefore,
$$
\int_{\B}\tilde{u}\Delta^{2}\varphi=\int_{\B}\frac{(1-\varepsilon)\mu}{\left(1-\frac{1}{1-\varepsilon}\tilde{u}\right)}dx
\geq (1-\varepsilon)\mu \int_{\B}\frac{1}{1-\tilde{u}}dx,\ \ \varphi\in C^{4}(\bar{\B})\cap H_{0}^{2}(\B)
$$
thus $\tilde{u}$ is a weak super-solution of
$(1.1)_{(1-\varepsilon)\mu}$ and Proposition \ref{P2.1} yields a
weak solution $v$ of $(1.1)_{(1-\varepsilon)\mu}$ which satisfies
$$
0\leq v\leq \tilde{u}<u_{\mu}\leq 1
$$
and then classical by Remark \ref{R1.1}.\vskip0.1in

\begin{remark}
From this Lemma, we know that $\lambda^{*}=\lambda_{*}$, in what
follows,  we always denote by $\lambda_{*}$ the largest possible
value of $\lambda$ such that $(1.1)_{\lambda}$ has a solution,
unless otherwise stated.
\end{remark}

\begin{proposition}\label{P2.2}
Up to a subsequence, the convergence
$$
u^{*}:=\lim_{\lambda\nearrow \lambda_{*}}\underline{u}_{\lambda}(x)
$$
holds in $H_{0}^{2}(\B)$ and the extremal solution $u_{\lambda^{*}}$
satisfies
\begin{equation}\label{Equ6}
\int_{\B}\Delta u^{*}\Delta\varphi=\lambda_{*}\int_{\B}\frac{1}{(1-u^{*})},\ \ \varphi\in C_{0}^{\infty}(\B).
\end{equation}
In particular, the extremal solution is weakly stable and if
$\|u^{*}\|_{\infty}<1$ then $\mu_{1}(u^{*})=0$.
\end{proposition}


{\bf Proof.} Since $\underline{u}_{\lambda}$ is stable, we have
\begin{equation}\label{Equ7}
\lambda\int_{\B}\frac{\underline{u}_{\lambda}^{2}}{(1-\underline{u}_{\lambda})^{2}}dx\leq\int_{\B}|\Delta\underline{u}_{\lambda}|^{2}dx
=\int_{\B}\underline{u}_{\lambda}\Delta^{2}\underline{u}_{\lambda}=\lambda\int_{\B}\frac{\underline{u}_{\lambda}}{1-\underline{u}_{\lambda}}dx.
\end{equation}
Next, it is easy to check that the following elementary inequality
holds: there exists a constant $C>0$ such that
$$
(1+C)\frac{s}{(1-s)}\leq\frac{s^{2}}{(1-s)^{2}}+(1+C), \ \  s\in (0,1),
$$
which used in (\ref{Equ7}) yields
$$
\lambda\int_{\B}\frac{\underline{u}_{\lambda}}{1-\underline{u}_{\lambda}}\geq\lambda\int_{\B}
\frac{\underline{u}_{\lambda}^{2}}{(1-\underline{u}_{\lambda})^{2}}dx\geq\lambda(1+C)
\int_{\B}\frac{\underline{u}_{\lambda}}{1-\underline{u}_{\lambda}}-C_{1},
$$
where $C_{1}$ is independent of $\lambda$. From the above
inequality, we get
$$
\|\Delta\underline{u}_{\lambda}\|_{2}^{2}=\lambda\int_{\B}\frac{\underline{u}_{\lambda}}{1-\underline{u}_{\lambda}}dx\leq C.
$$
Therefore, we may assume $\underline{u}_{\lambda}\rightharpoonup
u^{*}$ in $H_{0}^{2}(\B)$ and by monotone convergence theorem
(\ref{Equ6}) holds after integration by parts. Since
$\mu_{1}(\underline{u}_{\lambda})>0$ for all $\lambda\in
(0,\lambda_{*})$, in particular we have
$$
\int_{\B}|\Delta \varphi|^{2}dx\geq \int_{\B}\frac{\lambda \varphi^{2}}{(1-\underline{u}_{\lambda})^{2}},
\varphi\in C_{0}^{\infty}(\B)
$$
and passing to the limit as $\lambda\nearrow \lambda_{*}$ we obtain
that $u_{\lambda_{*}}$ is weakly stable. Finally, if
$\|u_{\lambda_{*}}\|_{\infty}<1$ and hence $u_{\lambda_{*}}$ is a
classical solution of $(1.1)_{\lambda_{*}}$, the linearized operator
at $u_{\lambda_{*}}$
$$
L(\lambda_{*}, u_{\lambda_{*}}):=\Delta^{2}-\frac{\lambda_{*}}{(1-u_{\lambda_{*}})^{2}}
$$
well defined on the space $\R^{+}\times C^{4,\alpha}(\B)$. If
$\mu_{1}(u_{\lambda_{*}})>0$ then the Implicit Function Theorem
applied to the function
$$
F(\lambda, \underline{u}_{\lambda}):=\Delta^{2} \underline{u}_{\lambda}-\frac{\lambda}{1-\underline{u}_{\lambda}}
$$
would yield a solution for $\lambda>\lambda_{*}$ contradicting the
definition of $\lambda_{*}$, thus $\mu_{1}(u^{*})=0$. \vskip0.1in


\begin{corollary}\label{C2.1}
 There exists a constant $C$ independent of
$\lambda$ such that for each $\lambda\in (0,\lambda_{*})$, the
minimal solution $\underline{u}_{\lambda}$ satisfies
$\|(1-\underline{u}_{\lambda})^{-1}\|_{L^{2}}\leq C$.
\end{corollary}

{\bf Proof.} From Proposition \ref{P2.2}, we have
$$
\int_{\B}\frac{\underline{u}_{\lambda}^{2}}{(1-\underline{u}_{\lambda})^{2}}dx=\int_{\underline{u}_{\lambda}\geq\frac{1}{2}}
\frac{\underline{u}_{\lambda}^{2}}{(1-\underline{u}_{\lambda})^{2}}dx
+\int_{\underline{u}_{\lambda}<\frac{1}{2}}\frac{\underline{u}_{\lambda}^{2}}{(1-\underline{u}_{\lambda})^{2}}dx\leq C.
$$
So
$$
\int_{\underline{u}_{\lambda}\geq\frac{1}{2}}
\frac{1}{(1-\underline{u}_{\lambda})^{2}}dx\leq4\int_{\underline{u}_{\lambda}\geq\frac{1}{2}}
\frac{\underline{u}_{\lambda}^{2}}{(1-\underline{u}_{\lambda})^{2}}dx\leq C.
$$
From this, we easily obtain
$\|(1-\underline{u}_{\lambda})^{-1}\|_{L^{2}}\leq C$, and the proof
is completed. \vskip0.2in

\begin{corollary}\label{C2.2}
 For dimensions $n\leq 4$, the
extremal solution $u^{*}$ is regular, i.e.,
$u^{*}=\lim_{\lambda\nearrow\lambda_{*}}u_{\lambda}$ exists in the
topology of $C^{4}(\B)$.
\end{corollary}

{\bf Proof } Since $u^{*}$ is radial and radially decreasing, we
need just to show that $u^{*}(0)<1$ to get the regularity of
$u^{*}$. Since $(1-u^{*}(x))\in L^{2}(\B)$  according to the
corollary \ref{C2.1}, we have that $u^{*}(x)\in W^{4,2}(\B)$ by the
standard elliptic regularity theory. And then by the Sobolev
imbedding theorem we have $u^{*}(x)\in C^{4-[\frac{n}{8}]-1,
[\frac{n}{8}]+1-\frac{n}{8} }(\B)$. So if $n\leq 4$, one can easy to
see that $u^{*}(x)\in C^{2}(\B)$. As $\nabla u^{*}(0)=0$, we get
$$
1-u^{*}(x)=u^{*}(0)-u^{*}(x)\leq C|x|^{2},
$$
hence
$$
\infty>\int_{\B}\frac{dx}{(1-u^{*}(x))^{2}}\geq C\int_{\B}\frac{dx}{|x|^{4}}=\infty.
$$
A contradiction arises, so $u^{*}$ is regular for $n\leq 4$.  \qed

\vskip0.2in
\subsection{  The upper and lower bounds for
$\lambda_{*}$ }\vskip0.2in

\begin{lemma}\label{L2.2}
$$
\lambda_{*}\leq \frac{\nu_{1}}{4},
$$
where $\nu_{1}$ is the first eigenvalue of $\Delta^{2}$ in
$H_{0}^{2}(\B)$. \vskip 0.1in
\end{lemma}
{\bf Proof.} Let $\underline{u}_{\lambda}$ be a solution of
$(1.1)_{\lambda}$ and let $(\psi,\nu_{1})$ denote the first
eigenpair of $\Delta^{2}$ in $H_{0}^{2}(\B)$ with $\psi>0$ then,
$$
\int_{\B}\underline{u}_{\lambda}\psi dx=\int_{\B}\underline{u}_{\lambda} \Delta^{2}\psi dx=\lambda\int_{\B}\frac{\psi}{1-\underline{u}_{\lambda}}
$$
and this implies
$$
\int_{\B}(-\nu_{1}\underline{u}_{\lambda}+\frac{\lambda}{1-\underline{u}_{\lambda}})\psi dx=0.
$$
Since $\psi>0$ there must exists a point $\bar{x}\in \B$ where
$$
\frac{\lambda}{1-\underline{u}_{\lambda}}-\nu_{1}\underline{u}_{\lambda}\leq 0,
$$
one can conclude that $\lambda_{*}\leq \sup_{0\leq
\underline{u}_{\lambda}\leq1}\nu_{1}\underline{u}_{\lambda}(1-\underline{u}_{\lambda})=\frac{\nu_{1}}{4}$.
\vskip0.1in

The lower bound for $\lambda_{*}$ is obtained by finding a suitable
supersolution . For example, if for some parameter
$\tilde{\lambda}_{1}$ there exists a supersolution, then
$\lambda_{*}>\tilde{\lambda}_{1}$ by Proposition \ref{P2.1}. \qed
\vskip 0.1in

\begin{lemma}\label{L2.3}
For $n\geq1$, we have
$$
\lambda_{*}\geq \max\{4n(n-2), 2n(n+2)\}.
$$
\end{lemma}

{\bf Proof.} For any $\beta>0$ and $C_{0}>0$ let
$g_{\beta}(r)=(C_{0}-\log r)^{\beta}, \ r\in (0,1)$. Then, by direct
calculation we find the following facts:
$$
\Delta g_{\beta}(r)=\beta r^{-2}[(2-n)g_{\beta-1}+(\beta-1)g_{\beta-2}];
$$
and
$$
\Delta[r^{2} g_{\beta}]=2n g_{\beta}-\beta(n+2)g_{\beta-1}+\beta(\beta-1)g_{\beta-2}.
$$
So we have
\begin{eqnarray*}
\Delta^{2}(r^{2}g_{\beta})&=&2n\Delta g_{\beta}-\beta(n+2)\Delta g_{\beta-1}(r)+\beta(\beta-1)\Delta g_{\beta-2}(r)\\\\
&=&\beta r^{-2} \times\left\{2n(2-n)g_{\beta-1}(r)+(\beta-1)(2n+n^{2}-4)g_{\beta-2}(r)\right\}\\\\
&+&\beta r^{-2} \times\left\{(\beta-1)(\beta-2)\times(-2n)g_{\beta-3}(r)+(\beta-1)(\beta-2)(\beta-3)g_{\beta-4}\right\}.
\end{eqnarray*}
Now let $\beta\in (0,1)$ and $n>2$, we have
\begin{equation} \label{Equ8}
\Delta^{2}(r^{2}g_{\beta})\leq \beta r^{-2}\times 2n(n-2)g_{\beta-1}.
\end{equation}

Also for any $A>0$ take $\bar{u}=1-Ar^{2}g_{\beta}$, one conclude
from (\ref{Equ8}) that
$$
\Delta^{2}u\geq 2n(n-2)A\beta r^{-2}g_{\beta-1}.
$$
Set $\beta=\frac{1}{2}$,  one can obtain that
$$\left\{
\begin{array}{lllllll}
\Delta^{2}\bar{u}\geq\frac{n(n-2)A^{2}}{1-\bar{u}} & \mbox{in} \ \  \B_{1}, \\\\
\bar{u}(r)=1-C_{0}^{\frac{1}{2}}A             & \mbox{on}\ \  \partial \B_{1},  \\\\
 \bar{u}'(r) =AC_{0}^{-\frac{1}{2}}(\frac{1}{2}-2C_{0} )   & \mbox{on}\ \ \partial \B_{1}.
 \end{array}
\right.
$$

Choosing $C_{0}=\frac{1}{4}, A_{0}=2$, one conclude that
$\bar{u}(r)$ is a supersolution of $(1.1)_{4n(n-2)}$ and
$\lambda_{*}\geq 4n(n-2)$ according to Proposition \ref{P2.1}.
Besides, we consider the function
$$
\omega_{\alpha}(x):=\alpha(1-|x|^{2})^{2},\ \  \alpha\in (0,1),
$$
which satisfies $0\leq \omega_{\alpha}(x)<1$ for $x\in \B$ and
$$
\omega_{\alpha}(x)=0, \frac{\partial\omega_{\alpha}}{\partial
n}=0,\  \mbox{for}\ \ x\in \partial \B;\ \mbox{for all}\ \ \alpha\in (0,1).
$$
Now the idea is to obtain from $\omega_{\alpha}(x)$ a super-solution
of $(1.1)_{\lambda}$, for a suitable choice of $\alpha$ and for
$\lambda$ in a suitable range of the form $0<\lambda\leq
\tilde{\lambda}$. For simply calculation, we have
\begin{eqnarray*}
\Delta^{2}\omega_{\alpha}(r)&=&\frac{d^{4}\omega_{\alpha}}{dr^{4}}+\frac{2(n-1)}{r}+\frac{d^{3}\omega_{\alpha}}{dr^{3}}
+\frac{(n-)(n-3)}{r^{2}}\frac{d^{2}\omega_{\alpha}}{dr^{2}}-\frac{(n-1)(n-3)}{r^{3}}\frac{d\omega_{\alpha}}{dr}\\
&=&[8n^{2}+16n]\alpha=:C(n)\alpha,
\end{eqnarray*}
and thus
$$
\Delta^{2}\omega_{\alpha}(r)=\frac{C(n)\alpha (1-\alpha)}{1-\alpha}\geq \frac{C(n)\alpha (1-\alpha)}{[1-\alpha(1-|x|^{4})]}
=\frac{C(n)\alpha (1-\alpha)}{1-\omega_{\alpha}}
$$
from which we deduce that
$$
\lambda_{*}=\lambda^{*}\geq \sup_{\alpha\in (0,1)}C(n)\alpha(1-\alpha)=\frac{1}{4}C(n)=2n(n+2)
$$
and the proof is completed.\qed \vskip0.1in

We complete this section by giving proofs of Theorem \ref{result1}
and \ref{result2}. \vskip0.1in

{\bf Proofs of Theorem \ref{result1} and \ref{result2}.}  The proof
of Theorem \ref{result1} follows form Proposition \ref{P2.1}, Lemma
\ref{L2.2} and Lemma \ref{L2.3}. For the proof of Theorem
\ref{result2},  we only need to prove the uniqueness of the regular
extremal solution $u^{*}$, the other parts of Theorem 1.2 follow
from Lemma \ref{L2.1} and Corollary \ref{C2.2}. Indeed, if the
extremal solution $u^{*}$ is regular, we can easily check that
$\mu_{1}(u^{*})=0$ by Implicit Function Theorem, since otherwise, we
can continue the minimal branch beyond $\lambda_{*}$. And then the
uniqueness follows from the (ii) of the Lemma \ref{L1.5}.
\qed\vskip0.2in

\setcounter{equation}{0}
 \section{Uniqueness of the extremal solution: proof of Theorem \ref{result3}}

\quad \quad {\bf Proof of Theorem \ref{result3}.} Suppose that $v\in
H^{2}(\B)$ satisfies
$$
\left\{
\begin{array}{lllllll}
\int_{\B} v\Delta^{2}\varphi dx\geq\int_{\B} \frac{\lambda_{*}}{1-v}dx, \forall \varphi\in C_{0}^{\infty}(\bar{\B}), \varphi\geq0,\\\\
v|_{\partial\B}=0,  \frac{\partial v}{\partial n} |_{\partial\B} \leq0,
\end{array}
\right.
$$
and $v\not\equiv u^{*}$. Notice that the construction of minimal
solutions in Proposition \ref{P2.1} for $\lambda\in
(0,\lambda_{*})$, carries over to $\lambda=\lambda_{*}$ but just in
the weak sense; precisely, we may assume that for
$\lambda=\lambda^{*}$ there exists a minimal weak solution. In other
words, it is legitimate to assume
$$
v(x)\geq u^{*},\ \  a.e.\ \ \ x\in \B.
$$
The idea of the proof is as follows: first we prove the function
$$
u_{0}=\frac{1}{2}(u^{*}+v)
$$
is a super-solution to the following perturbation of problem
$(1.1)_{\lambda}$
\begin{equation}\label{Equ10}
\left\{
\begin{array}{lllllll}
 \Delta^{2} u=\frac{\lambda_{*}}{1-u}+\mu\frac{\zeta(x)}{1-u}, & \mbox{in}\ \ \B;\\
 0\leq u\leq1, & \mbox{in}  \ \ \B;\\
 u=\frac{\partial u}{\partial n}=0, & \mbox{on}  \ \ \partial\B,\\
 \end{array}
\right.
\end{equation}
for a standard cut-off function $\zeta(x)\in C_{0}^{\infty}(\B)$ and
$\mu>0$ to be suitably chosen;  besides, a solution is understood in
weak sense unless otherwise stated. Second, we construct, for some
$\lambda>\lambda_{*}$, a super-solution to $(1.1)_{\lambda}$ by
using a solution of (\ref{Equ10}) and this will enable us to build
up a weak solution of $(1.1)_{\lambda}$ for $\lambda>\lambda_{*}$
and thus necessarily $v\equiv u^{*}$.

Indeed we first observe that for $0<R<1$ and for some
$c_{0}=c_{0}(R)>0$
\begin{equation}\label{Equ11}
v(x)\geq u^{*}+c_{0}\ \ \ |x|\leq R.
\end{equation}
To prove this we recall the Green's function for $\Delta^{2}$ with
Dirichlet boundary conditions
$$
\left\{
\begin{array}{lllllll}
 \Delta_{x}^{2} G(x,y)=\delta_{y}&\ \ \ x\in\B; \\
 G(x,y)=0  &\ \ \ x\in\partial\B; \\
 \frac{\partial G}{\partial n} (x,y)=0 &\ \ \ x\in\partial\B,
 \end{array}
\right.
$$
where $\delta_{y}$ is the Dirac mass at $y\in \B$. Boggio gave an
explicit formula for $G(x,y)$ which was used in \cite{Gru} to prove
that in dimension $n\geq 5$
\begin{equation}\label{Equ11}
G(x,y)\sim|x-y|^{4-n}\min(1, \frac{d(x)^{2}d(y)^{2}}{|x-y|^{4}}),
\end{equation}
where
$$
d(x)=\mbox{dist}(x,\partial\B)=1-|x|.
$$
Formula (\ref{Equ11}) yields
\begin{equation}\label{Equ12}
G(x,y)\geq cd(x)^{2}d(y)^{2}
\end{equation}
for some $c>0$ and this in turn implies that for smooth functions
$\bar{v}$ and $\bar{u}$ such that $\bar{v}-\bar{u}\in H_{0}^{2}(\B)$
and $\Delta^{2}(\bar{v}-\bar{u})\geq0$,
\begin{eqnarray*}
\tilde{v}-\tilde{u}&=&\int_{\partial B}\left(\frac{\partial\Delta_{x}G}{\partial n_{x}}(x,y)\tilde{v}-\tilde{u}
-\Delta_{x}G(x,y)\frac{\partial (\tilde{v}-\tilde{u})}{\partial n}\right)\\
&+&\int_{B}G(x,y)\Delta^{2}(\tilde{v}-\tilde{u})dx\\
&\geq&c d(y)^{2}\int_{B}\Delta^{2}(\tilde{v}-\tilde{u}) d(x)^{2}dx.
\end{eqnarray*}
Using a standard approximation procedure, we conclude that
$$
v(y)-u^{*}(y)\geq cd(y)^{2}\lambda^{*}\int_{\B}\left(\frac{1}{1-v}-\frac{1}{1-u^{*}}\right)d(x)^{2}dx.
$$
Since $v\geq u^{*}, v\not\equiv u^{*}$ we deduce (\ref{Equ11}).

Let $u_{0}=\frac{u^{*}+v}{2}$. Then by Taylar's Theorem
\begin{equation}\label{Equ13}
\frac{1}{1-v}=\frac{1}{1-u_{0}}+\frac{v-u_{0}}{(1-u_{0})^{2}}+\frac{(v-u_{0})^{2}}{4(1-u_{0})^{3}}
+\frac{(v-u_{0})^{3}}{18(1-u_{0})^{3}}+\frac{(v-u_{0})^{4}}{96(1-\varepsilon_{1})^{4}}
\end{equation}
for some $u_{0}\leq \varepsilon_{1}\leq v$ and
\begin{equation}\label{Equ14}
\frac{1}{1-u^{*}}=\frac{1}{1-u_{0}}+\frac{u^{*}-u_{0}}{(1-u_{0})^{2}}+\frac{(u^{*}-u_{0})^{2}}{4(1-u_{0})^{3}}
+\frac{(u^{*}-u_{0})^{3}}{18(1-u_{0})^{3}}+\frac{(u^{*}-u_{0})^{4}}{96(1-\varepsilon_{2})^{4}}
\end{equation}
for some $u^{*}\leq\varepsilon_{2}\leq u_{0}$. Adding (\ref{Equ13})
and (\ref{Equ14}) yields
\begin{equation}\label{Equ15}
\frac{1}{2}(\frac{1}{1-v}+\frac{1}{1-u^{*}})\geq\frac{1}{1-u_{0}}+\frac{1}{16}\frac{(u^{*}-v)^{2}}{(1-u_{0})^{2}}
\end{equation}
and in turn we obtain,
\begin{eqnarray*}
\int_{\B}u_{0}\Delta^{2}\varphi dx\geq \int_{\B}\left[\frac{\lambda_{*}}{1-u_{0}}+\frac{\lambda_{*}(u^{*}-v)^{2}}{16(1-u_{0})}\right]dx
\geq\int_{\B}\left[\frac{\lambda_{*}}{1-u_{0}}+\frac{\lambda_{*}c_{0}^{2}\zeta(x)}{16(1-u_{0})}\right]dx.
\end{eqnarray*}
Thus, $u_{0}$ is a weak super-solution of (\ref{Equ10}) with
$\mu=\frac{\lambda_{*}c_{0}^{2}}{16}$ and the cut-off $\zeta(x)$
with support in $\B_{\rho}$. Now reasoning as in Lemma \ref{L2.1},
we may assume for $\varepsilon>0$ sufficiently small, that
(\ref{Equ10}) posses a classical solution $0\leq u_{\varepsilon}<1$
with parameter $\lambda_{*}$ replaced by $\lambda_{*}-\varepsilon$.
Set $\mu_{\varepsilon}:=[(\lambda_{*}-\varepsilon)c_{0}^{2}]/16$ and
let $\psi\in C^{4}(\bar{\B})$ be the unique classical solution of
the following
$$
\left\{
\begin{array}{lllllll}
\Delta^{2}\psi=\mu_{\varepsilon} \frac{\zeta(x)}{1-u_{\varepsilon}} &\mbox{in}\ \ \B,\\
\psi=\frac{\partial \psi}{\partial n}=0 &\mbox{on}\ \ \partial\B.
\end{array}
\right.
$$
We also, by the Boggio principle, have that there exists $M>0$
sufficiently large such that $u_{\varepsilon}\leq M\psi$. Next let
$\delta>0$ and set
$$
\omega:=\frac{(\lambda_{*}-\varepsilon)+\delta}{\lambda_{*}-\varepsilon}u_{\varepsilon}-\psi
$$
and choosing $\delta$ sufficiently small, we obtain $\omega\leq
u_{\varepsilon}<1$; moreover, from
$$
\left\{
\begin{array}{lllllll}
\Delta^{2}(u_{\varepsilon}-\psi)=(\lambda_{*}-\varepsilon)\frac{1}{1-u_{\varepsilon}}\geq0,&\ \ \mbox{in}\ \ \B,\\
 u_{\varepsilon}-\psi=\frac{\partial(u_{\varepsilon}-\psi)}{\partial n}=0  &\ \ \mbox{on} \ \ \partial\B,
 \end{array}
\right.
$$
we have again by the Boggio principle that $\psi\leq
u_{\varepsilon}$ and eventually that $\omega\geq 0$. Finally we have
$$
\Delta^{2}\omega=(\lambda_{*}-\varepsilon+\delta)\frac{1}{1-u_{\varepsilon}}
+\frac{(\lambda_{*}-\varepsilon+\delta)c_{0}^{2}}{16}\frac{\varepsilon(x)}{1-u_{\varepsilon}}
-\mu_{\varepsilon}\frac{\varepsilon(x)}{1-u_{\varepsilon}}\geq (\lambda_{*}-\varepsilon+\delta)\frac{1}{1-\omega}
$$
since $\omega\leq u_{\varepsilon}$. Thus it is enough to choose
$0<\varepsilon<\delta$ to provide a classical solution to
$(1.1)_{\lambda}$ for $\lambda>\lambda_{*}$ which is a
contradiction; this completes the proof of Theorem \ref{result3}.

 \vskip0.2in
\setcounter{equation}{0}
\section{Behavior of the minimal solutions as
$\lambda\to0$: proof of Theorem \ref{result4}}\vskip0.1in

\quad \quad { \bf Proof of  Theorem \ref{result4}.}\ \  We first
show that
\begin{equation}\label{Eq4.1}
\underline{u}_{\lambda}\to 0\ \ \mbox{uniformly as }\ \ \lambda\to 0. 
\end{equation}
Since this standard, we just briefly sketch its proof. By Theorem
\ref{result1}, we know that
$$
0<\lambda<\mu<\lambda_{*}\Rightarrow \underline{u}_{\lambda}(x)<\underline{u}_{\mu}(x)
\ \ \mbox{if}\ \ |x|< 1.
$$
Then, by multiplying the equation $(1.1)_{\lambda}$ by
$\underline{u}_{\lambda}$ and by integrating by parts, we obtain
that $\|\underline{u}_{\lambda}\|_{H_{0}^{2}(\B)}$ remains bounded.
Hence, up to a subsequence, $\{\underline{u}_{\lambda}\}$ converges
in the weak $H_{0}^{2}(\B)$ topology to 0, which is the unique
solution of $(1.1)_{0}$. By convergence of the norms, we infer that
the convergence is in the norm topology.

Next, note that $U_{\lambda}$ satisfies
$$
\left\{
\begin{array}{lllllll}
\Delta^{2}U_{\lambda}=\lambda & \mbox{in}\ \B, \\
U_{\lambda}=\frac{\partial U_{\lambda}}{\partial n}  =0 &   \mbox{on}\ \partial\B.
 \end{array}
\right.
$$
Therefore,
$\Delta^{2}\underline{u}_{\lambda}>\Delta^{2}U_{\lambda}$,  one
conclude that $\underline{u}_{\lambda}>U_{\lambda}$ by Lemma 2.1.

In order to prove the last statement of Theorem \ref{result4}, note
that from (\ref{Eq4.1}) we know that
$$
\mbox{for all}\ \varepsilon>0\ \ \mbox{there exists}\ \lambda_{\varepsilon}>0\ \ \mbox{such that}
\lambda<\lambda_{\varepsilon}\Rightarrow \|\underline{u}_{\lambda}\|_{\infty}<\varepsilon.
$$
So, fix $\varepsilon>0$ and let $\lambda<\lambda_{\varepsilon}$.
Then
$$
\Delta^{2}\underline{u}_{\lambda}=\frac{\lambda}{1-\underline{u}_{\lambda}}<\frac{\lambda}{1-\varepsilon}
=\Delta^{2}\frac{U_{\lambda}}{1-\varepsilon}
$$
This shows that
$\underline{u}_{\lambda}(x)<\frac{U_{\lambda}(x)}{1-\varepsilon}$
for all $x\in \B$, and the proof is completed according to the
arbitrariness of $\varepsilon$. \qed

\setcounter{equation}{0}
\section{Further results and open problems}
\quad \quad First, we give the following result which is the main
tool to guarantee that $u^{*}$ is singular. At the same time, it
give a precise estimate for $\lambda_{*}$. The proof of this result
is based on an upper estimate of $u^{*}$ by a stable singular
subsolution.

\begin{proposition}\label{P4.1}
Suppose there exist $\lambda'>0, \beta>0$ and a singular radial
function $\omega(r)\in H_{0}^{2}(\B)$ with $\frac{1}{1-\omega(r)}\in
L_{loc}^{\infty}(\bar{\B}\setminus 0)$ such that
\begin{equation}\label{Eq5.1}
\left\{
\begin{array}{lllllll}
 \Delta^{2} \omega\leq\frac{\lambda'}{1-\omega}& \mbox{for}\ \ 0<r<1,\\
  \omega(1)=  \omega'(1)=0,
 \end{array}
\right.
\end{equation}
and
\begin{equation}\label{Eq5.2}
\beta\int_{\B}\frac{\phi^{2}}{(1-\omega)^{2}}\leq \int_{\B}(\Delta\phi)^{2}\quad \mbox{for all}\quad \phi\in H_{0}^{2}(\B).
\end{equation}
If $\beta>\lambda'$, then $\lambda_{*}<\lambda'$ and $u^{*}$ is
singular.
\end{proposition}

{\bf Proof.}\quad First, note that (\ref{Eq5.2}) and
$\frac{1}{1-\omega(r)}\in L_{loc}^{\infty}(\bar{\B}\setminus 0)$
yield to $\frac{1}{1-\omega}\in L^{1}(\B)$. (\ref{Eq5.1}) implies
that $\omega(r)$ is a $H_{0}^{2}(\B)-$ weak sub-solution of
$(1.1)_{\lambda'}$. If now $\lambda'<\lambda^{*}$, then by Lemma
\ref{L1.5}, $\omega(r)$ would necessarily be below the minimal
solution $\underline{u}_{\lambda'}$, which is a contradiction since
$\omega(r)$ is singular while $\underline{u}_{\lambda'}$ is regular.
In the following, we shall prove that $u^{*}$ is singular.

Now let $\frac{\lambda'}{\beta}<\gamma<1$ in such a way that
$$
\alpha:=(\frac{\gamma \lambda_{*}}{\lambda'})^{\frac{1}{2}}<1.
$$
Setting $\bar{\omega}:=1-\alpha(1-\omega)$, we claim that
\begin{equation}\label{Eq5.3}
u^{*}\leq \bar{\omega} \ \ \ \mbox{in}\ \ \ \B. 
\end{equation}
Note that by the choice of $\alpha$ we have
$\alpha^{2}\lambda'<\lambda_{*}$, and therefore to prove
(\ref{Eq5.3}) it suffices to show that for $\alpha^{2}\lambda'\leq
\lambda<\lambda_{*}$, we have $u_{\lambda}\leq \bar{\omega}$ in
$\B$. Indeed, fix such $\lambda$ and note that
$$
\Delta^{2}\bar{\omega}=\alpha \Delta^{2}\omega\leq\frac{\alpha \lambda'}{(1-\omega)}
=\frac{\alpha^{2}\lambda'}{(1-\bar{\omega})}\leq \frac{\lambda}{(1-\bar{\omega})}.
$$
Assume that $\underline{u}_{\lambda}\leq \bar{\omega}$ dose not hold
in $\B$, and consider
$$
R_{1}:=\sup\{0\leq R\leq 1| \underline{u}_{\lambda}(R)>\bar{\omega}(R)\}>0.
$$
Since $\bar{\omega}(1)=1-\alpha>0=u_{\lambda}(1)$, we then have
$$
R_{1}<1, u_{\lambda}(R_{1})=\bar{\omega}(R_{1})\  \mbox{and}\  \underline{u}_{\lambda}'(R_{1})\leq\bar{\omega}'(R_{1}).
$$
Now consider the following problem
$$
\left\{
\begin{array}{lllllll}
\Delta^{2}u=\frac{\lambda}{1-u}& \ \mbox{in}\ \ \B_{R_{1}},\\
u=u_{\lambda}(R_{1})& \ \  \mbox{on}\ \ \partial \B_{R_{1}},\\
 \frac{\partial u}{\partial n}= u'_{\lambda}(R_{1}) & \ \  \mbox{on}\ \ \partial \B_{R_{1}}.
\end{array}
\right.
$$
Then $\underline{u}_{\lambda}$ is a solution to above problem while
$\bar{\omega}$ is a sub-solution to the same problem. Moreover
$\bar{\omega}$ is stable since $\lambda<\lambda_{*}$ and
$$
\frac{\lambda}{(1-\bar{\omega})^{2}}\leq\frac{\lambda_{*}}{\alpha^{2}(1-\omega)^{2}}
<\frac{\beta}{(1-\omega)^{2}}.
$$
By Lemma \ref{L2.1}, we deduce that $\underline{u}_{\lambda}\geq
\bar{\omega}$ in $\B_{R_{1}}$ which is impossible, since
$\bar{\omega}$ is singular while $u_{\lambda}$ is regular. This
establishes claim (\ref{Eq5.3}) which, combined with the above
inequality, yields
$$
\frac{\lambda_{*}}{(1-u_{*})^{2}}\leq\frac{\lambda_{*}}{\alpha^{2}(1-\omega)^{2}}
<\frac{\beta}{(1-\omega)},
$$
and thus
$$
\inf_{\varphi\in C_{0}^{\infty}(\B)}\frac{\int_{\B}[(\Delta \varphi)^{2}-\frac{\lambda_{*}\varphi^{2}}{(1-u^{*})^{2}}]dx}{\int_{\B}\varphi^{2}dx}>0.
$$
This is not possible if $u^{*}$ is a smooth function, since
otherwise, one could use the Implicit function Theorem
 to continue the minimal branch beyond $\lambda_{*}$. The proof is
 over. \qed

\begin{itemize}
\item $Open\ Problem\ 1$. Dose $(1.1)_{\lambda}$ exist a stable singular
subsolution? We know that Cowan etal, with the help of  Maple,
construct such solution of $(P_{\lambda})$ with $p=2$ by
improved Improved Hardy-Rellich Inequalities, see \cite{Mo,Co1}
. But the method used there seems invalid.
\end{itemize}

We now turn to the extremal solution $u^{*}$. We  suggest the
following open problems.
\begin{itemize}
\item $Open\ Problem\ 2$. Dose one find the precise estimate for $u^{*}$ as
in \cite{Mo,Co1,Da}, which  play a crucial role for
investigating the regularity of $u^{*}$. In  \cite{Co1}, the
precise bound for $u^{*}$ is obtained by finding a stable
singular subsolution which relies on the \textquotedblleft
ghost" singular solution, as mentioned in introduction. However,
in the present paper we can not find any \textquotedblleft
ghost" singular solution, so a new trick is needed.
\end{itemize}

\begin{itemize}
\item $Open\ Problem\ 3$. For the corresponding second equation the
extremal solution $u^{*}$ is regular for dimensions $n\leq6$ and
singular for dimension $n\geq7$, for details see \cite{Mea}. The
threshold $n^{*}=7$ between regular and singular solutions is
called the critical dimension. There is a natural question:
whether there exists a critical dimension $N^{*}$ for equation
$(1.1)_{\lambda}$. We conjecture that $N^{*}=8$. \qed

\end{itemize}
\noindent \textbf{Acknowledgement.} The author is greatly indebted
to Prof. Dong Ye for his constructive comments. This research is
supported in part  by National Natural Science Foundation of China
(Grant No. 10971061).


\small {\it 1 Institute of Contemporary Mathematics, Henan University;}\\
\small {\it 2  School of Mathematics and Information Science,Henan University}\\
\small {\it Kaifeng 475004, P. R. China.}


\begin{thebibliography}{99}

\bibitem{Pe} {\small J. A. Pelesko, A. A. Bernstein, Modeling MEMS and NEMS, Chapman Hall and CRC Press,
2002.}
\bibitem{Es} {\small P. Esposito, Compactness of a nonlinear eigenvalue problem with a singular nonlinearity,
Commun. Contemp. Math. 10 (2008), no. 1, 17-45}
\bibitem{Es1} {\small  P. Esposito, N. Ghoussoub and Y. Guo, Compactness along the branch of semi-stable and unstable
solutions for an elliptic problem with a singular nonlinearity,
Comm. Pure Appl. Math. 60 (2007),  1731-1768.}
\bibitem{Gh} {\small N. Ghoussoub, Y. Guo, On the partial differential equations of electro MEMS devices: stationary case,
SIAM J. Math. Anal. 38 (2007), 1423-1449.}

\bibitem{Ca} {\small D. Cassani, J. do O ,  N. Ghoussoub, On a fourth order elliptic problem with a singular nonlinearity, Advances Nonlinear
Studies, 9, (2007), 177-197.}

\bibitem{Co1} {\small C. Cown, P. Esposito, N. Ghoussoub, and A. Moradifam, The critical dimension for a forth
order elliptic problem with singular nonlineartiy, Arch. Ration.
Mech. Anal., (2010, to appear).}

 \bibitem{D} {\small J. D\`{a}vila, I. Flores, I. Guerra, Multiplicity of solutions for a fourth order equation with power-type nonlinearity,
 Math. Ann. 348, (2009), 143-193.}

 \bibitem{Gu4} {\small Z. Guo, J. Wei, On a fourth order nonlinear elliptic equation with negative exponent, SLAM J. Math. Anal. 40, (2009),
2034-2054.}

\bibitem{Li} {\small F. Lin , Y. Yang, Nonlinear non-local elliptic equation modelling electrostatic actuation, Proc. R. Soc. A 463,
(2007), 1323-1337.}
\bibitem{Mo} {\small  Amir. Moradifam,  On the critical dimension of a fourth order elliptic problem with negative exponent,
Journal of Differential Equations 248 (2010), 594-616.}

\bibitem{Bo} {\small T. Boggio, Sulle funzioni di Freen d¡¯ordine m. Rend. Circ. Mat. Palermo 20 (1905), 97-135.}
\bibitem{Ghe} {\small M. Ghergu, A biharmonic equation with singular nonlinearity, http://arxiv.org/abs/0911.0308. }

\bibitem{Co1} {\small C. Cown, P. Esposito, N. Ghoussoub, and A. Moradifam, The critical dimension for a forth order elliptic problem with singular
nonlineartiy, Arch. Ration. Mech. Anal. (2009, to appear).}

\bibitem{Mea} {\small  A. M. Meadows, Stable and Singular Solutions of the Equation $\Delta u=\frac{1}{u}$,  Indiana Univ Math. J (53), 1681-1703.}

\bibitem{Ag} {\small S. Agmon, A. Dougist, and L. Nirenberg, Estimates near the boundary for solutions of elliptic partial differential
equations satisfying general boundary conditions. I, Comm. Pure
Appl. Math. (12), 1959, 623¨C727.}

\bibitem{Br} {\small H. Brezis, T. Cazenave, Y. Martel, A. Ramiandrisoa, Blow up for $u_{t}-\Delta u=g(u)$ revisited, adv.
Differential Eqations, (I), 1996, 73-90.}

\bibitem{Ma} {\small Y. Martel, Uniqueness of weak extremal solutions for nonlinear elliptic problems, Houston J. Math., 23, 161-168 (1997).}
\bibitem{Da} {\small J. D\`{a}vila, L. Dupaigne, I. Guerra, and M. Montenegro, Stable Solutions for the Bilaplacian with Exponential nonlinearity,
Siam J. Math. Anal., 39 (2007), 565-592.}
\bibitem{Ar1} {\small G. Arioli, F. Gazzola, H.-C. Grunau, E. Mitidieri, A semilinear fourth order elliptic problem
with exponential nonlinearity. Siam J. Math. Anal. 36, (2005),1226-1258.}

\bibitem{Gru} {\small H. Ch. Grunau and G. Sweers, Positivity for equations involving polyharmonic operators
with Dirichlet boundary conditions, Math. Ann. (307) 1997, 589¨C626.
}




\bibitem{Mi} {\small  F. Mignot, J. P. Puel, Solution radiale singuli\`{e}re de $-\Delta u=\lambda e^{u}$, C. R. Acad. Sci.
 Paris S\'{e}r. I 307, (1988), 379-382.}
\bibitem{Br1} {\small Brezis, H., Vazquez, J. L., Blow up solutions of some nonlinear elliptic problems.
 Rev. Mat. Univ. Complutense Madrid 10, (1997), 443-469.}





\end{thebibliography}
\end{document}